\documentclass[12pt]{article}

\usepackage{amsmath,amsfonts,amssymb,amsthm ,mathrsfs,setspace}

\newcommand{\exn}{{\bf E}}%{\mbox{\boldmath$E$}}
\newcommand{\pr}{{\bf P}}%{\mbox{\boldmath$P$}}
\newcommand{\N}{\mathbb{N}}
\newcommand{\R}{\mathbb{R}}
\newcommand{\Z}{\mathbb{Z}}

\newcommand{\Bi}{\mathscr{B}}
\newcommand{\Ii}{\mathscr{I}}
\newcommand{\Ai}{\mathcal{A}}

\newcommand{\Vi}{\mathcal{V}}

\newcommand{\Si}{\mathcal{S}}

\newcommand{\al}{\alpha}

\newcommand{\si}{\sigma}

\newcommand{\deq}{\stackrel{d}{=}}
\newcommand{\ep}{\varepsilon}

\newcommand{\ind}{{\bf 1}}

\newcommand{\bX}{\mbox{\boldmath$X$}}     % bold X
\newcommand{\bY}{\mbox{\boldmath$Y$}}     % bold Y
\newcommand{\bZ}{\mbox{\boldmath$Z$}}     % bold Z\
\newcommand{\bmu}{\mbox{\boldmath$\mu$}}     % bold mu
\newcommand{\bnu}{\mbox{\boldmath$\nu$}}     % bold nu

\newcommand{\bla}{\mbox{\boldmath$\lambda$}}     % bold la
\newcommand{\bka}{\mbox{\boldmath$\varkappa$}}     % bold kappa
\newcommand{\beeta}{\mbox{\boldmath$\eta$}}     % bold  eta
\newcommand{\bx}{\mbox{\boldmath$x$}}

\newcommand{\bm}{\mbox{\boldmath$m$}}
\newcommand{\bt}{\mbox{\boldmath$t$}}
\newcommand{\bs}{\mbox{\boldmath$s$}}
\newcommand{\bp}{\mbox{\boldmath$p$}}
\newcommand{\bu}{\mbox{\boldmath$u$}}     % bold u\

\newcommand{\sbX}{\mbox{\scriptsize\boldmath$X$}}     % small bold X
\newcommand{\sbt}{\mbox{\scriptsize\boldmath$t$}}     % small bold t
\newcommand{\sbnu}{\mbox{\scriptsize\boldmath$\nu$}}     % small bold nu
\newcommand{\sbla}{\mbox{\scriptsize\boldmath$\lambda$}}   % small bold lambda

\newtheorem{thm}{Theorem}
\newtheorem{coro}{Corollary}

\newtheorem{rema}{Remark}

\begin{document}

\title{On the asymptotic behaviour of a dynamic version of the Neyman contagious point process%
\footnote{Research supported by the ARC Discovery Grant DP120102398.}
}

\author{K.~Borovkov\footnote{Department of Mathematics and Statistics, The University of Melbourne, Parkville 3010, Australia; e-mail: {borovkov@unimelb.edu.au}.}
}

\date{}

\maketitle

\begin{abstract}
We consider a dynamic version of the Neyman contagious point process that can be used
for modelling the spacial dynamics of biological populations, including species
invasion scenarios. Starting with an arbitrary finite initial configuration of points
in $\R^d$ with nonnegative weights, at each time step a point is chosen at random from
the process according to the distribution with probabilities proportional to the
points' weights. Then a finite random number of new points is added to the process, each
displaced from the location of the chosen ``mother'' point by a random vector and assigned a
random weight. Under broad conditions on the sequences of the numbers of newly added
points, their weights and displacement
vectors (which include a random environments setup), we derive the asymptotic behaviour of the locations of the points added to the process at time
step~$n$ and also that of the scaled mean measure of the point process after time
step~$n\to\infty$.

\smallskip

{\em Key words and phrases:} Neyman contagious point process, random network, preferential attachment,
random environments, stationary sequence, limit theorems, law of large numbers, the
central limit theorem, regular variation.

\smallskip
{\em AMS Subject Classification 2010:} 60G55 % Point processes
 (Primary), 60F05 % Central limit and other weak theorems
 (Secondary).

\end{abstract}

\section{Introduction and main results}

The present paper deals with a dynamic  version of the Neyman contagious point process~\cite{Ne39, Th54} (the term ``contagious" in the context of such point processes going back to G.~P\'olya \cite{Po31}). The latter can be described as follows. Suppose we are given a homogeneous Poisson point process on the carrier space $\R^d,$ $d\ge 1$ (in a  motivating example from~\cite{Ne39}, the points representing egg masses of some insect species). Then, using each of the points in the process as a ``mother'', we add several ``daughter" points randomly displaced relative to their mother (to model a situation where ``larvae hatched from the eggs which are being laid in so-called `masses'"~\cite{Ne39}). Our process differs from that one-stage scheme in that it is multistage and, at each time step, one of the existing points is chosen at random to be the next mother, and then a random number of daughter points are added, displaced at random relative to their mother point. The points in our process are endowed with ``weights" that determine the probabilities with which the points can be chosen to be mothers in the future.

More formally, our point process starts with an initial configuration of $k_0\ge 1$  random points
$\bX_{0,1},\ldots,\bX_{0,k_0}\in \R^d,$ $d\ge 1,$  labelled with respective vectors
$(w_{0,j}, u_{0,j}),$ $j=1, \ldots,  k_0,$ with non-negative components.

At time step   $n\ge 1,$ $k_n\ge 0$ new points are added to the process, with respective labels
$(w_{n,j}, u_{n,j}),$ $j=1, \ldots,  k_n.$ The role of the weights $w_{n,l}\ge 0$ will
be to characterise the ``fitness", or reproductive ability, of the individuals to be associated
with the points to be added to the process at the $n$th step (they can be used, for
example, to model the aging of the individuals), as they will be used to determine the
probability distributions for choosing new mothers in the future. The quantities
$u_{n,l}\ge 0$ specify the amount of a ``resource" attached to the individuals (e.g.,
the body weight etc.).

To specify the locations of the newly added points, let $  W_{-1}:=0,$ and, for $n\ge 0$, set
\[
w_n:=   \sum_{l=1}^{k_n} w_{n,l}, \qquad W_n:= \sum_{r=0}^n w_r,
\]
Then, at step $n+1\ge 1$, an existing point $\bX_{r,j},$ $0\le r\le n,$ $1\le j \le k_r$, is chosen at random,  according  the probability distribution
\begin{equation}
\bp_{n }:=\{ p_{n,r,j} : 0\le r\le n, \, 1\le j \le k_r \},\qquad p_{n,r,j}:=\frac{w_{r,j}}{W_{n }}.      \label{ps}
\end{equation}
Denote that point by $\bX^*_n$ and add $k_{n+1} \ge 0$ new points  to the process, at the locations
\[
 \bX_{n+1, l}:=\bX^*_n + \bY_{n+1,l},\qquad   1\le l\le k_{n+1},
\]
where we assume that the random vector $(\bY_{n+1,1}, \ldots, \bY_{n+1,k_{n+1}})\in(\R^d)^{k_{n+1}}$ is independent of $\{\bX_{r,j}\}_{0\le r \le n,1\le j\le k_r }$ and the choice of the ``mother point" $\bX_n^*$. We assume nothing about the character of dependence between the components $\bY_{n+1,l}$ (in particular, some of them can coincide).
Denote the distribution of $(\bY_{n,1}, \ldots, \bY_{n,k_n})$ on $(\R^d)^{k_n}$ by ${\bf Q}_n$, and that  of $\bY_{n,j}$ on $\R^d$ by~$Q_{n,j}$.

Note that if the $w_n$'s tend to decrease as $n$ increases, it means that the more mature individuals are more active in reproduction. When the $w_n$'s tend to increase, the younger ones are more productive.

Models of such kind can be used to describe the dynamics of a population of microorganisms in varying (in both time and space) environments, and the process of distribution of  individual insects (as in the  larvae example from~\cite{Ne39}) or colonies of social insects (such as some bee species or ants). In particular, they can be of interest when modeling foreign species invasion scenarios. Further examples include dynamics of business development for companies employing  multi-level marketing techniques (or referral marketing).

Observe that our process' structure   resembles that of a branching random walk, see,
e.g.,~\cite{AsKa76, Bi97} and references therein. If $\bp_n$ is a  uniform distribution, the process under consideration is  the embedded skeleton process for a continuous time branching random walk. There is also an interesting relationship with certain random search algorithms, see, e.g., \cite{De88} and~\eqref{probab_interp} below.

One can think about the points in our process as nodes in a growing random forest, the roots thereof being the initial points, with edges connecting daughters to their mothers. Note that the trees in the forest, unlike the ones in the case of branching processes, are not independent of each other.

The model we are dealing with belongs to the class of ``preferential attachment processes". The most famous of them was proposed in the much-cited paper~\cite{BaAl99}, where new nodes were attached to already existing ones with probabilities proportional to the degrees of the latter, with the aim of   generating random ``scale-free" networks that  are widely observed in natural and human-made systems, including the World Wide Web and some social networks. Since then there have  appeared quite a few publications on the topic, including the monograph~\cite{Du07}. In most cases, in the existing literature the ``attachment rule" depends on the degrees of the existing nodes (possibly with some variations as, e.g., in~\cite{Aietal09}, where  new nodes may link to a
node only if they fall within the node's ``influence region"), and the authors consider the  standard set of questions concerning the growing random graph, such as: the appearance and size of the giant component (if any), the sizes of (other) connected components in different regimes, the diameter and average distance between two vertices in the giant component,  typical degree distributions, proportion of vertices with a given degree, the maximum degree, bond percolation and critical probability, and connectivity properties after malicious deletion of a certain proportion of vertices.

In relation to our previous remark on connection to branching phenomena, we would like to mention the paper~\cite{RuTyVa07}   obtaining the asymptotic degree distribution in a preferential attachment random tree for a wide range of weight functions (of the nodes' degrees) using well-established results from the theory of general branching processes.

However, to the best of the author's knowledge, neither models of the type discussed  in the present paper nor the questions of the spacial dynamics of the growing networks had been considered in the literature prior to the publication of~\cite{BoMo05}. That paper dealt with the special case of the model where $k_n\equiv k=\mbox{const},$ $n\ge 1$, the weights being  assumed
to be of the form $w_{n,j}=a^n,$ $j=1,\ldots,k ,$ for some constant~$a>0,$ with
$\bY_{n,j},$ $j=1,\ldots, k,$ being independent and identically distributed.  Under rather broad conditions
(existence of Ces\`aro limits for means/covariance matrices of $\bY_{n,1}$, uniform
integrability for $\bY_{n,1}$ or $\|\bY_{n,1}\|^2$, $\|\cdot\|$ denoting the Euclidean norm), the paper established versions  of
our Theorems~\ref{Thm1}--\ref{Thm3} below. They showed that the distributions of the
locations of the individuals added at  the $n$th step and the mean measures of the
process after step~$n$ display   distinct LLN- and CLT-type behaviours depending on whether $a<1,$ $a=1$ or
$a>1.$ The model was further studied in~\cite{BoVa06}, in the case where $k=1,$ $d=1,$
$\bY_{n,1}\ge 0,$ and the weights $w_n$ were assumed to be random. In addition to
obtaining  an analog of our Theorem~\ref{Thm1} and proving a number of other results,
the paper also established the asymptotic (as $n\to \infty$) behaviour  of the distance
$\bX_{n}^*+\bY_{n+1,1}$ to the root (the initial point in the process) of the point
added at step $n+1$ under the assumption that $w_n = e^{-S_n}$, $S_n=\theta_1+ \cdots +
\theta_n$, $\theta_j$ being i.i.d.\ with zero mean  and finite variance. The
paper~\cite{BoVa06} also derived the asymptotics of the expected values of the
outdegrees of vertices.

In the present work, we are dealing with a broader class of models, where the numbers $k_n$ of points added to the process at different steps can vary, while, at any given time step $n$, the displacement vectors $\bY_{n,j}$, $j=1,\ldots,   k_n$, can be dependent and have different distributions, and study not only the laws of the locations of the newly added points, but also the distributions of their ``resources" (that may coincide with their ``fertility" given by the weights~$w_{n,j}$).

Moreover, we will also be dealing with  processes of the above type evolving in random environments (RE). This  means that the offspring sizes $k_n$, their attributes $(w_{n,j}, u_{n,j})$ and the distributions of $(\bY_{n,1}, \ldots, \bY_{n,k_n})$ will be random as well, which can be formalised as follows.

Let
$\mathscr{P}$ be a given family of distributions on $( \R^d, \Bi (\R^d)),$ $ \Bi (\R^d)$ being the $\si$-algebra of Borel sets on~$ \R^d$. For $k\ge 1,$ denote by
$\mathscr{P}^{(k)}$ the family of all distributions on $(\R^d)^k$ with   margins  from $\mathscr{P}$.

The random environments are given    by a sequence $\Vi:= \{V_n\}_{n\ge 1}$ of
random elements $V_n  \in S$  of the space
\[
S:= \bigcup_{k\ge 0} S^{(k)}, \quad S^{(k)}:= (\R_+^2)^k\times \mathscr{P}^{(k)}, \quad
k=1,2,\ldots, \quad  S^{(0)}:=\{0\}.
\]
The elements $V_n$  can be specified by their ``dimensionality"  $k_n$ (one has $k_n=k$ when $V_n\in
S^{(k)}$, meaning that $k$ new points will be added to the process at step~$n$),
``attributes'' $((w_{n,1},u_{n,1}), \ldots, (w_{n,k_n},u_{n,k_n}) )\in (\R_+^2)^{k_n}$,
and offspring displacement distributions ${\bf Q}_n\in \mathscr{P}^{(k_n)}.$ The choice of the $\sigma$-algebra on $S$ is standard for objects of such
nature, so we will not describe it in detail.

The dynamics of the process in RE are basically the same as in the case of the deterministic environments, the only difference being that, in the above-presented description of step $n+1$ in the process, one should everywhere add ``given the whole sequence of random environments~$\Vi$ and the past history of the
process up to step $n $". Thus, $\bp_{n }$ will become the (random) conditional distribution (given~$\Vi$) of choosing a new mother among the already existing points $\bX_{0,1},  \ldots, \bX_{n,k_n}$ etc.

The conditions of our assertions in the RE setup  will often include strict stationarity of some sequences related to~$\Vi$. In those cases, $\Ii$ will be used to denote the $\si$-algebra of   events invariant w.r.t.\ the respective measure preserving transformation,
\[
\langle \,\cdot\, \rangle :=\exn  (\,\cdot \,| \,\Ii )
\]
denoting the conditional expectation given~$\Ii$.

To formulate the main results of the paper, we will need some further notation.  Denote by \[
\phi_{n,j} (\bt) := \exn \exp\{ i\bt\bX_{n,j}^T \},\quad \bt\in\R^d,
\]
the characteristic function (ch.f.) of the distribution $P_{n,j}$ of the  point $\bX_{n,j},$ $n\ge 0,$ $1\le j \le k_n,$ in our process, $ \bx^T$ standing for the transposed
(column) vector  of $\bx\in\R^d$, and set, for  $ n=0,1,2,\ldots$ and $\bt\in
\R^d,$
\begin{align}
\pi_n := \frac{w_n}{W_n}\equiv  \sum_{j=1}^{k_n} p_{n,n,j} ,
 \quad \phi_n (\bt) := \frac1{w_n} \sum_{j=1}^{k_n} w_{n,j}  \phi_{n,j}(\bt)\equiv \frac1{\pi_n} \sum_{j=1}^{k_n} p_{n,n,j}  \phi_{n,j}(\bt)
 \label{fpi}
\end{align}
(note that $\pi_0=1)$.

That is, $\phi_n$ is the ch.f.\ of the mixture law $  P_n:= \frac1{w_n} \sum_{j=1}^{k_n} w_{n,j}  P_{n,j}$, which coincides with the conditional distribution  of the location $\bX_n^*$ of the mother point for step $n+1$ in our process given that that point was chosen from the $k_n$ points added at the previous step. For $n$ such that $k_n=0,$ we can leave $\phi_{n,j}$ and $ \phi_n$ undefined.

In the RE setup, we put
\begin{align*}
\phi_{n,j|\Vi} (\bt) := \exn_{\Vi}  \exp\{i\bt \bX_{n,j}^T\}, \quad
\phi_{n|\Vi} (\bt) := \frac1{w_n} \sum_{j=1}^{k_n} w_{n,j}  \phi_{n,j|\Vi} (\bt)  ,\quad  \bt\in
\R^d,
\end{align*}
where  $\exn_{\Vi} $ stands for the conditional expectation given $\Vi $, and denote by $P_{n,j|\Vi}$ and  $P_{n|\Vi}$ the respective distributions.

Set $ u_n:=   \sum_{j=1}^{k_n} u_{n,j},$  $U_n:= \sum_{r=0}^n u_r,$ and consider the
measure
\begin{align}
\label{mean_m} M_n (B) &
 := \exn  \biggl[\frac{1}{U_n}\sum_{r=0}^n \sum_{j=1}^{k_r} u_{r,j}
  \ind_{\{\sbX_{ r,j} \in B\}}\biggr]
 = \frac{1}{U_n}  \sum_{r=0}^n \sum_{j=1}^{k_r} u_{r,j}  P_{r,j} (B), \quad B\in\Bi (\R^d),
\end{align}
describing the distribution of the ``resource" specified by the quantities~$u_{n,j}$,
$\ind_A$ being the indicator of the event~$A$. When $u_{n,j}\equiv 1,$ $M_n$ is just the mean measure of the process. In the RE setup, $M_{n|\Vi}$ is defined by~\eqref{mean_m} with $\exn$ in its middle part replaced  with~$\exn_{\Vi}$ (and hence with $P_{r,j}$ on its right-hand side replaced with $P_{r,j|\Vi}$).

To avoid making repetitive trivial comments, we assume throughout the paper that
$U_n\to \infty$ (a.s.\ in the RE setup) as $n\to\infty,$  so that the process of adding new points and resources never terminates.

Further, we will   denote by
\[
f_{n,j} (\bt) := \exn  \exp\{i\bt \bY_{n,j}^T\}, \quad \bt \in\R^d,
\]
the  ch.f.\ of $Q_{n,j}$ ($f_{n,j|\Vi} (\bt)$ is defined in the same way as $f_{n,j} (\bt)$, but with $\exn$ replaced with~$\exn_{\Vi}$), and set
\[
f_n (\bt) := \frac1{w_n} \sum_{j=1}^{k_n} w_{n,j}  f_{n,j}(\bt) .
\]

Denote by $\Ai$ the directed set of all pairs $(n,j),$ $n\ge 0, $ $1\le j\le k_n,$ with preorder $\preccurlyeq$ on it defined by $(r,j) \preccurlyeq (n,l)$ iff $r\le n$. Given a net $a:=\{a_\alpha\}_{\alpha\in\Ai}$ in a linear topological  space, and a net $\{b_\alpha\}_{\alpha\in\Ai}$ in $[0,\infty)$,
we say that the net $a$ is $b$-summable with sum $\Si_b (a)$ if the limit
\[
 \Si_b (a):= \lim_\alpha\frac{ \sum_{\beta\preccurlyeq \alpha}b_\beta a_\beta }{\sum_{\beta\preccurlyeq \alpha}b_\beta}
\]
exists and is finite.
If $b_{n,j}\equiv 1$ then $b$-summability is basically Ces\`aro summability (for nets), to a sum denoted by $\Si_1 (a)$. In the above definition,   $a$ can be a net of scalars, vectors, matrices, or even distributions. In the latter case the limit (to be denoted by $w\mbox{-}\lim$) is   in  the topology of weak convergence of distributions,   whereas in other cases the limit is in the component-wise (point-wise for functions) sense.

\begin{rema}\label{rem_sum}
{\rm An elementary extension of the classical Toeplitz theorem shows that, if $\lim_\alpha a_\alpha =a_\infty$ and $\lim_\alpha b_\alpha  =\infty$, then $a$ is $b$-summable with the sum $a_\infty$. Observe also that
if both $ab=\{a_\alpha b_\alpha\}$ and $b$ are 1-summable, and
$  \Si_1 (b) >0$, then $ a$ will be $b $-summable as well, with
the sum $\Si_b (a) = \Si_1 (ab)/\Si_1 (b).$}
\end{rema}

As we already said, the weights $w_{n,j}$ can be used to model dependence of the
reproductive ability of individuals on when there were added to process (in particular,
their ``aging"). Our first result refers to the case when the total ``weight'' of all
the points in the process is finite, meaning that the
more mature individuals tend to have noticeably higher reproduction ability compared to the newly added ones.

\begin{thm}
\label{Thm1} Let $W_\infty:=\lim_{n\to\infty} W_n\equiv \sum_{r=0}^\infty w_n <\infty$. Then  the infinite product
\[
\Pi (\bt) :=\prod_{r=0}^\infty (1 + \pi_r (f_r (\bt ) -1))
\]
converges  to a ch.f. Moreover,

{\rm (i)}~if $\{Q_{n,j}\}$ is $u$-summable $($or, which is the same, the net
$\{f_{n,j} (\bt)\}$ is $u$-summable to a continuous sum $ \Si_u ( f (\bt) ))$,  then one has  $w\mbox{-}\lim_{n\to\infty} M_n=  M_\infty$, where the limiting measure $M_\infty$ has ch.f.\
   $\Pi (\bt)  \Si_u ( f (\bt) );$

{\rm (ii)}~if $w\mbox{-}\lim_{(n,j)}  Q_{n,j}= Q_\infty $   for some distribution~$Q_\infty$ on $\R^d$ with ch.f.\ $f_\infty (\bt)$, then   one has $w\mbox{-}\lim_{(n,j)}  P_{n,j}= P_\infty$,  the limiting distribution $P_\infty$ having the ch.f.\
\[
\phi_\infty (\bt):=   f_\infty (\bt) \Pi (\bt).
\]
In that case, one also has  $w\mbox{-}\lim_{n\to\infty} M_n=  P_\infty$
\end{thm}

From the Birkhoff--Khinchin theorem, Remark~\ref{rem_sum} and   Theorem~\ref{Thm1}(i), we immediately obtain the following result in the RE setup. We write $\bZ \sim P$ if the the random vector $\bZ$ has distribution~$P$,   denote by  $\bu_n:=(u_{n,1}, \ldots, u_{n,k_n})\in \R_+^{k_n}$ the vector of resource values $u_{n,j}$ assigned to the points added to the process at step~$n$, and set
\[
(uf)_n (\bt):= \sum_{j=1}^{k_n} u_{n,j}f_{n,j} (\bt).
\]

\begin{coro}
In the RE setup, suppose that $\{\|\bZ\|:\bZ \sim P\in\mathscr{P} \}$ is uniformly integrable and $\{(k_n, \bu_n, {\bf Q}_n (\cdot))\}$ is a strictly stationary sequence.
Then, on the event $\{W_\infty <\infty\}\cap\{ \langle u_1\rangle >0\},$ one has $w\mbox{-}\lim_{n\to\infty} M_{n|\Vi}= M_\infty$ a.s., where the limiting measure $M_\infty$ has ch.f.
\begin{equation*}
    \Pi (\bt)\frac{ \langle (uf)_1 (\bt)\rangle}{\langle u_1\rangle}.
\end{equation*}
\end{coro}

An RE version of the assertion of Theorem~\ref{Thm1}(ii) concerning convergence of $P_{n|\Vi}$ is straightforward.

\begin{rema}
{\rm
Thus, when $W_\infty< \infty$ and the displacement laws $Q_{n,j}$ ``stabilise on average" in the sense that they are $1$-summable, the points in our process form a ``cloud" of which the relative ``density" has a limit. This is so because the ``old points" will be responsible for most of the progeny. But what about the contribution of individual points? How often will the points be chosen to become mothers (recall that we do not exclude the case $k_n=0$, so a chosen point will not necessarily have daughters at a given step)? Since, for any fixed $r\ge 0$ and $j=1,\ldots, k_r$, the indicators of the events  $\{\bX_{r,j}$ is chosen to be the mother at step $n+1\}$,  $n\ge r,$   are independent Bernoulli random variables with success probabilities $w_{r,j}/W_n$,  Borel--Cantelli's lemma immediately yields the following dichotomy:
\begin{itemize}
\item if $\sum_{n\ge 0}W_n^{-1}= \infty$ then any individual with a positive $w$-weight will a.s.\ be chosen to become a mother infinitely often (so this is the case in Theorem~\ref{Thm1}), and
\item if $\sum_{n\ge 0}  W_n^{-1} < \infty$ then any individual   will become a mother finitely often a.s.
\end{itemize}
}
\end{rema}

The latter situation is, of course, impossible under conditions of Theorem~\ref{Thm1}, but both situations can occur under the conditions of
the next theorem that deals with the case where the total weight of the points in the process is unbounded. For $n\ge 0$, set $\bmu_{n,j}:=\exn  \bY_{n,j},$ $\bm_{n,j}:=\exn  \bY_{n,j}^T\bY_{n,j} $ (when these expectations are finite),  and let
\[
 \bmu_n:=\frac1{w_n} \sum_{j=1}^{k_n} w_{n,j}\bmu_{n,j},
\quad
 \bm_n:=\frac1{w_n} \sum_{j=1}^{k_n} w_{n,j}\bm_{n,j}.
\]
It is clear that $\bmu_n$ is the mean vector of the distribution with ch.f.~$f_n,$ while $\bm_n$ is the matrix of (mixed) second order moments of that mixture distribution.

\begin{thm}
\label{Thm2} Assume that $w_n = \xi_n n^\alpha L(n)$, $n\ge 1,$  where $\{\xi_n\ge 0\}$ is Ces\`aro summable to a positive value $\Si_1 (\xi)$,  $\alpha>-1$ and $L$ is a slowly varying at infinity function. %, both being deterministic.

{\rm (i)} Let   $\{\|\bY_{n,j}\| \}_{(n,j)\in\Ai}$ be
uniformly integrable,   $\{\xi_n \bmu_n\}$ be Ces\`aro summable with sum $\Si_1 (\xi \bmu)$.   Then
\begin{equation}
 w\mbox{-}\!\lim_{(n,j)} P_{n,j} (\, \cdot\, \ln n)  = \delta_{\sbla}(\, \cdot\, ),
 \label{thm2_1}
\end{equation}
where $\delta_{\sbla}$  is the unit mass concentrated at the point
$\bla:=(\alpha+1) {\Si_1 (\xi \bmu) }/{\Si_1 (\xi)}\in\R^d,$
and
\[
\bnu_n:=\frac1{\ln n}\sum_{r=1}^{n-1} \pi_r  \bmu_r \to \bla, \quad n\to\infty.
\]
Moreover, if
\begin{equation}
\lim_{\ep\searrow 0}\limsup_{n\to\infty} U_{\ep n}/U_n=0 ,
\label{Uep}
\end{equation}
then also
\begin{equation}
 w\mbox{-}\!\!\lim_{n\to\infty}M_n (\, \cdot\, \ln n)  = \delta_{\sbla}(\, \cdot\, ).
  \label{thm2_2}
\end{equation}

{\rm (ii)}  Assume that $\{\|\bY_{n,j}\|^2 \}_{(n,j)\in\Ai}$ is uniformly integrable and
the sequence   $\{\xi_n \bm_n\}$  is Ces\`aro summable. Then the limit $w\mbox{-}\!\!\lim_{(n,j)}$ of the  distributions  of
\[
\frac{\bX_{n,j} - \bnu_n \ln n}{\sqrt{\ln n}}
\]
exists and is equal to the normal law in $\R^d$ with zero mean vector
and covariance matrix $(\alpha+1)  {\Si_1 (\xi \bm)}/{\Si_1 (\xi )}.$
Moreover, if \eqref{Uep} holds then $w\mbox{-}\!\!\lim_{n\to\infty} M_n (\, \cdot\,\sqrt{\ln n} +\bnu_n\ln n)$ also exists and is equal to the same normal distribution.
\end{thm}

As in the case of Theorem~\ref{Thm1}, an RE version of the above statement is readily available from  the Birkhoff--Khinchin theorem and Remark~\ref{rem_sum}. Set
\[
\bmu_{n|\Vi}:=\frac1{w_n}\sum_{j=1}^{k_n} w_{n,j}\exn_{\Vi} \bY_{n,j},
\qquad \bm_{n|\Vi} :=\frac1{w_n}\sum_{j=1}^{k_n} w_{n,j}\exn_{\Vi} \bY_{r,j}^T \bY_{r,j}.
\]

\begin{coro}
In the RE setup, let $w_n = \xi_n n^\alpha L(n)$, $n\ge 1,$ where $\alpha=\alpha (\omega)>-1$ and $L(n) = L(\omega,n)$ is a slowly varying function a.s.

{\rm (i)} Assume that $\{\|\bZ\|: \bZ\sim P \in \mathscr{P}\}$ is
uniformly integrable and
\begin{equation}
\label{cond_xi}
\xi_n= \widetilde\xi_n (1 + o(1)) , \quad
\bmu_{n|\Vi}= \widetilde\bmu_n  + o(1) \quad\mbox{a.s.\ as $n\to\infty$,}
\end{equation}
where
$\{(\widetilde\xi_n, \widetilde\bmu_n)\}$ is  a strictly stationary sequence with a finite absolute moment and  $\exn \widetilde\xi_1   \|\widetilde\bmu_1\| <\infty$. Then, as $n\to\infty,$ on the event $A:= \{\langle\widetilde \xi_1\rangle >0\},$ one has
\begin{equation*}
 w\mbox{-}\!\lim_{(n,j)} P_{n,j|\Vi} (\, \cdot\, \ln n)  = \delta_{\sbla}(\, \cdot\, ) \quad
 \mbox{a.s.,}
\end{equation*}
where $\delta_{\sbla}$  is the unit mass concentrated at the point
$\bla:=(\alpha+1) {\langle \widetilde\xi_1 \widetilde \bmu_1 \rangle}/{\langle\widetilde\xi_1\rangle}\in\R^d$, and
\[
\bnu_{n|\Vi}:=\frac1{\ln n}\sum_{r=1}^{n-1} \pi_r \bmu_{r|\Vi} \to \bla\quad \mbox{a.s.}
\]
Moreover, if \eqref{Uep} holds a.s.\ on $A$, then on that event one also has
\begin{equation*}
 w\mbox{-}\!\!\lim_{n\to\infty}M_{n|\Vi}  (\, \cdot\, \ln n)  = \delta_{\sbla}(\, \cdot\, )\quad
 \mbox{a.s.}
\end{equation*}

{\rm (ii)}  Assume that the family $\{\|\bZ\|^2: \bZ\sim P \in \mathscr{P}\}$ is uniformly integrable and the first of relations \eqref{cond_xi} holds together with
\begin{equation}
 \label{m_n}
\bm_{n|\Vi} =\widetilde \bm_n + o(1)\quad\mbox{\it  a.s.\ \ as \ } n\to\infty,
\end{equation}
where
the sequence  $\{(\widetilde\xi_{n},   \widetilde\bm_n)\}$  is strictly stationary with  a finite  absolute moment and $\exn \widetilde\xi_1 \|\widetilde \bm_1\|<\infty$.  Then, on the event $A,$   the limit $w\mbox{-}\!\!\lim_{(n,j)}$ of the  conditional distributions  of
\[
\frac{\bX_{n,j} - \bnu_{n|\Vi} \ln n}{\sqrt{\ln n}}
\]
given $\Vi $ exists  a.s.\ and equals the normal law in $\R^d$ with zero mean vector
and covariance matrix $(\alpha+1) {\langle \widetilde\xi_1 \widetilde \bm_1 \rangle}/{ \langle \widetilde\xi_1\rangle}.$ Moreover, if \eqref{Uep} holds on $A$,  that normal distribution will be limiting on $A$  for $  M_{n|\Vi} (\, \cdot\,\sqrt{\ln n} +\bnu_n\ln n)$  as well.
\end{coro}

Thus, in the broad (power function) range of the $w$-weights dynamics, the point processes display  basically one and the same behaviour: a drift at rate $\ln n$ in the direction of the ``average" mean displacement vector (if the latter is non-zero), and the asymptotic normality of the point distribution at the scale $\sqrt{\ln n}$. The next theorem deals with the case where the $w$-weights grow exponentially fast on average, basically meaning that the individuals modelled by  points in our process  have a ``finite life span" (at least, in what concerns their reproduction abilities).  As one could expect in such a situation, the distributions $P_n$ will ``drift" at a constant rater, while the $u$-resource will be ``spread" along the trajectory of the drift.

In this case, the formulation of the result is much simpler (and more natural) in the RE setup than in the deterministic one, so we will only  consider here the former situation.

\begin{thm}
\label{Thm3}
In the RE setup, let $w_n = \xi_n e^{S_n}$, $n\ge 1,$ where $S_n:=\tau_1+ \cdots + \tau_n$ for a random sequence $\{\tau_n\}_{n\in\Z}$.

{\rm (i)} Assume that $\{\|\bZ\| : \bZ\sim P \in \mathscr{P}\}$ is
uniformly integrable,  relations   \eqref{cond_xi} holds with  $ \{( \widetilde\xi_n, \tau_n,\widetilde\bmu_n)\}_{n\in\Z}$ being   strictly stationary, $\exn \widetilde\xi_1 \| \widetilde\bmu_1\|<\infty.$ Then, on the event $A:= \{\langle \widetilde\xi_1\rangle >0, \langle \tau_1\rangle >0\},$
\begin{equation}
 \label{Pn}
w\mbox{-}\!\lim_{(n,j)} P_{n,j|\Vi} (\, \cdot\, \ln n)  = \delta_{\sbla}\quad \mbox{a.s., \quad}
 \bla := \Bigl\langle \frac{ \widetilde\xi_1 \widetilde\bmu_1}{ \widetilde\zeta_1}\Bigr\rangle,
\end{equation}
where
\[
\widetilde\zeta_n:=\sum_{m=0}^\infty  \widetilde\xi_{n-m}e^{-S_{n,m}},
 \quad
S_{n,m}:= \tau_{n-m+1}+\tau_{n-m+2} + \cdots + \tau_n,\quad m>1,
\]
$S_{n,0}:=0,$ so that $ \{( \widetilde\xi_n, \widetilde\zeta_n,\widetilde\bmu_n)\}$   is   stationary   as well on $A$.

If, in addition, the distribution functions $  U_{ \lfloor nv\rfloor   }/U_n, $ $v\in [0,1]$, converge weakly a.s.\  to  some distribution function $G$ on $[0,1]$ as $n\to\infty,$  then,  on the event $A$ there exists the limit  $w\mbox{-}\!\lim_{n\to\infty} M_{n|\Vi} (\, \cdot\, n)= M_\infty (\, \cdot\,)$  a.s., where $M_\infty$ is supported by the straight line segment with end points ${\bf 0}$ and $\bla,$ such  that~$M_\infty (\{\bla s, s\in [0,v]\})=G(v),$    $v\in [0,1].$

{\rm (ii)}  Let    $\{\|\bZ\|^2 : \bZ\sim P \in \mathscr{P}\}$   be uniformly integrable,   relations \eqref{cond_xi} and  \eqref{m_n} hold   with the sequence  $\{( \widetilde\xi_{n}, \tau_n, \widetilde\bmu_n,  \widetilde\bm_n)\}_{n\in\Z}$ being strictly stationary with a finite
first absolute moment.  Then, on the event $A,$   the limit $w\mbox{-}\!\!\lim_{(n,j)}$ of the  conditional given $\Vi $ distributions  of
\[
\frac{\bX_{n,j} - \bka_n  }{\sqrt{  n}}, \qquad    \bka_n:= \sum_{r=1}^n \pi_r\bmu_{r|\Vi},
\]
exists a.s.\ and equals  the normal law in $\R^d$ with zero mean vector and covariance matrix $\Bigl\langle \dfrac{  \widetilde\xi_1 }{ \widetilde\zeta_1}\, \widetilde\bm_1 -  \dfrac{  \widetilde\xi_1^2  }{ \widetilde\zeta_1^2}\,\widetilde\bmu^T_1\widetilde\bmu_1\Bigr\rangle.$

\end{thm}

\begin{rema}{\rm
In part~(ii) of Theorem~\ref{Thm3}, one can also derive a normal mixture approximation to $M_{n|\Vi} (\, \cdot\, \sqrt{n})$ as $n\to\infty$, cf. Theorem~1(iii) in~\cite{BoMo05}.}
\end{rema}

\section{Proofs}

The key tool in the proofs is an extended version of the recurrences (9)
in~\cite{BoMo05} and~(2) in~\cite{BoVa06}. Namely, first note that,  by the total probability formula,
\begin{align}
\label{phi2}
\phi_{n+1,j} = \Phi_n f_{n+1,j},
 \quad \Phi_n :=  \sum_{r=0}^n \sum_{l=1}^{k_r}p_{n , r,l} \phi_{r,l}.
\end{align}
Using~\eqref{ps} and the observation that $p_{n,r , j} =p_{n-1,r, j}  W_{n-1}/W_{n}=p_{n-1,r, j}(1-\pi_n) $ for all $r<n$, $1\le j\le k_r,$ we obtain for $n\ge 1 $,   employing~\eqref{phi2}, \eqref{fpi} and recursion, that
\begin{align}
\Phi_n &= \sum_{r=0}^{n-1} \sum_{l=1}^{k_r} p_{n , r,l} \phi_{r,l}
 +\sum_{l=1}^{k_n} p_{n , n,l} \phi_{n,l}
\notag\\
 &=  (1-\pi_n) \Phi_{n-1} +   \sum_{l=1}^{k_n} p_{n , n,l}\Phi_{n-1} f_{n,l}
 = (1-\pi_n) \Phi_{n-1} + \pi_n \Phi_{n-1} f_n
 \notag\\&
  = (1+\pi_n (f_n -1) ) \Phi_{n-1}
   = \prod_{r=0}^n (1+\pi_r (f_r -1) ).
   \label{forma}
\end{align}
A probabilistic interpretation of \eqref{forma} is that
\begin{equation}
\bX_{n}^*\deq \sum_{r=0}^n I_r \bY_r  ,\quad n\ge 1,
\label{probab_interp}
\end{equation}
where the $I_r$'s are Bernoulli random variables with success probabilities $\pi_r$ (so that $I_0\equiv 1$), the random vectors $\bY_r $ have ch.f.'s $f_r$, all of the random quantities being independent.

Therefore, it follows from \eqref{phi2} that
\begin{align}
\phi_{n+1,j}   =    f_{n+1,j}  \prod_{r=0}^n \bigl(  1 + \pi_r (f_r  -1)\bigr),\quad n\ge 0.
\label{phi3}
\end{align}
In the RE setup,   fixing  $\Vi $ and using the same argument, we obtain that
$\phi_{n+1,j|\Vi}$ also equals the right-hand side of the above relation.

\begin{rema}
{\rm Representation~\eqref{probab_interp} basically corresponds to ``going back in time", tracing the ``ancestry" of the points added at step~$n+1$.
In the case where $w_n\equiv 1$, $\bY_n\equiv 1$, this is equivalent to the correspondence between the ``distance to the root" in a random tree and the number of records in an i.i.d.\ sequence of random variables used in~\cite{De88}. Observe also the following:  fix an arbitrary $(r,i)\in \Ai.$ Then, for any point $\bX_{n,j}$ with $n>r$,  the probability that $\bX_{n,j}$  has $\bX_{r,i}$ among its ancestors is one and the same value $w_{r,i}/W_r.$ Hence the mean number of the $n$th generation points that have a particle from the $r$th generation as an ancestor is $k_n\pi_r$.
}
\end{rema}

{\em Proof of Theorem~\ref{Thm1}.}  Since in this case $\pi_n = (1+ o(1))w_n/W_\infty$ as
$n\to\infty,$  it follows that $\sum_{r=1}^\infty \pi_r<\infty.$ Hence the product on
the right-hand side of~\eqref{forma} converges   to $\Pi (\bt) $ uniformly in $\bt$ as $n\to\infty$. That the limiting infinite product will be a ch.f.\ follows from L\'evy's continuity theorem, since all the factors $1 + \pi_r (f_r  -1)$ are ch.f.'s (note
that $\pi_r\in (0,1)$) and  the Weierstrass theorem implies that
$\Pi (\bt) $  will be continuous at zero.

(i)~To prove convergence of~$M_n$,  observe that~\eqref{mean_m}, the assumption that $U_n\to\infty$ as $n\to\infty$ and~\eqref{phi3} imply that, for any fixed $\bt\in\R^d,$ for the ch.f.\ $\varphi_n (\bt)$ of $M_n$ one has
\begin{align*}
\varphi_n (\bt) &
 =  \frac{1}{U_n}  \sum_{r=0}^n \sum_{j=1}^{k_r} u_{r,j}  \phi_{r,j} (\bt)
 =   \frac{1}{U_n}  \sum_{r=0}^n \sum_{j=1}^{k_r} u_{r,j} f_{r,j}  (\bt)
  \prod_{s=0}^{r-1} \bigl(  1 + \pi_s (f_s  (\bt)  -1)\bigr)\\
& = (1+o(1))\Pi (\bt) \frac{1}{U_n} \sum_{r=0}^n \sum_{j=1}^{k_r} u_{r,j} f_{r,j}(\bt)
 = (1+ o(1)) \Pi (\bt) \Si_u (f(\bt)).
\end{align*}
Now the desired assertion immediately follows from     L\'evy's continuity theorem.

(ii)~That $w\mbox{-}\lim_{(n,j)}  P_{n,j}= P_\infty$ follows from representation~\eqref{phi3},  the assumption that $ \lim_{(n,j)}  f_{n,j} (\bt)  = f_\infty (\bt)$ and the already established convergence of the products of $1 + \pi_r (f_r  -1)$  to the ch.f.\ $\Pi (\bt).$   That $w\mbox{-}\lim_{n\to\infty} M_n= P_\infty$ follows from the first part of Remark~\ref{rem_sum} and part~(i) of the theorem.\hfill$\Box$\medskip

{\em Proof of Theorem~\ref{Thm2}.} (i)~First we  will analyse the asymptotic behaviour
of $\pi_n$ as~$n\to\infty.$ Choose a sequence  $\varepsilon_p \searrow 0$ as $p\to\infty$ that vanishes slowly enough, so that the following relations hold true:  $m_p\to\infty$ as $p\to\infty$ for the sequence $\{m_p\}$ defined by   $m_{-1}:=-1$, $m_0:=1,$
\[
 m_p:= (1+\varepsilon_p) m_{p-1} ,\quad p=  1,2,\ldots,
\]
$\ep_p/\ep_{p-1}\to 1$, and
\begin{align}
\vartheta_p:= \frac1{m_p}\sum_{r=1}^{m_p}\xi_r - \Si_1 (\xi) = o(\ep_p).
  \label{theta0}
\end{align}
It is clear that such a sequence exists, and
we can assume for notational simplicity   that all $m_j$ are integer (the
changes needed in case they are not are obvious) and that there exists a $q=q(n)\in
\N$ such that $m_q=n$. Then
\begin{align*}
\sum_{r=1}^n w_r = \sum_{p=0}^{q }
  \sum_{r\in (m_{p-1},m_{p}]}  \xi_r r^{\alpha} L (r)
 = \sum_{p=0}^{q } m_p^{\alpha} L (m_p)
 \sum_{r\in (m_{p-1},m_{p}]} \xi_r \left(\frac{r}{m_p}\right)^\alpha \frac{L(r)}{L (m_p)},
\end{align*}
and it is not hard to verify (using the fact that $\sum_n n^\alpha L(n)=\infty$) that the right-hand side here (and hence $W_n$ as well) is
\begin{align}
(1+ o(1)) \sum_{p=0}^{q } m_p^{\alpha} L (m_p)
 \sum_{r\in (m_{p-1},m_{p}]}\xi_r.
\label{summa}
\end{align}
From Ces\`aro summability of $\{\xi_n\}$ and~\eqref{theta0}, one has
\begin{align*}
\frac1{\ep_p m_{p-1}}\sum_{r\in (m_{p-1},m_{p}]} \xi_r
 &= \frac1{\ep_p m_{p-1}}\biggl(\sum_{r=  1}^{m_{p}} \xi_r
  - \sum_{r=  1}^{m_{p-1}}\xi_r \biggr)
  \\
  &= \frac1{\ep_p } \bigl[(1+\ep_p)
  \bigl( \Si_1 (\xi) +\vartheta_{p}\bigr) - \bigl( \Si_1 (\xi) + \vartheta_{p-1}\bigr)\bigr]
  \\
  &
   =  \Si_1 (\xi) +\vartheta_{p} + (\vartheta_{p}-\vartheta_{p-1})/{\ep_p}
   =  \Si_1 (\xi) + o(1).
\end{align*}
Therefore  the representation~\eqref{summa} for $W_n$ implies that
\begin{align*}
W_n &= (1+ o(1))  \Si_1 (\xi) \sum_{p=1}^{q } \ep_p  m_p^{\alpha +1 } L (m_p)
     \\
    & = (1+ o(1))  \Si_1 (\xi)\sum_{r=1}^{n} r^{\alpha} L (r)
 = \frac{(1+ o(1))  \Si_1 (\xi) }{\alpha+1}\,   n^{\alpha+1} L (n)
\end{align*}
by Karamata's theorem, so that
\begin{align}
\pi_n \equiv  \frac{w_n}{W_n }
 = (\alpha+1 + o(1)) \frac{\xi_n}{ n  \Si_1 (\xi) }=o(1) \quad \mbox{   as $n\to\infty$},
 \label{beta_n}
\end{align}
where the last relation holds since $\xi_n=o(n)$  from Ces\`aro summability of $\{\xi_n\}.$

Now we   turn to   analysing the asymptotic behaviour of $P_{n,j}$. For a fixed
$\bt\in\R^d$ and a sequence $b_n\to\infty,$ $n\to\infty,$ it follows from~\eqref{phi3}
that the ch.f.\ of $P_{n,j} (\cdot~b_n)$  is given by
\begin{align}
\phi_{n,j} (\bt /b_n)
  & =   f_{n,j } (\bt/b_n)  \prod_{r=0}^{n-1} \bigl[  1 + \pi_r
(f_r (\bt/b_n) -1)\bigr]
 \notag\\
& =
 (1+o(1))  f_{n,j } (\bt/b_n) \exp\left\{  (1+o(1)) \sum_{r=0}^{n-1} \pi_r (f_r (\bt/b_n) -1)\right\}
 \notag\\
 &= (1+o(1)) \exp\left\{  \frac{1+o(1)}{b_n} \sum_{r=1}^{n-1} \pi_r (\bmu_r\bt^T
 + \eta_{r,n})  \right\},\quad |\eta_{r,n}|\le \eta_n=o(1),
\label{forma_b}
\end{align}
where the second equality follows from the relation
\[
\lim_{n\to\infty}\sup_{r \ge 0} |\pi_r (f_r (\bt/b_n) -1)| \to 0,
\]
which holds due to~\eqref{beta_n}, and the third one follows from the uniform
integrability assumption on~$\mathscr{P}$ (cf.~(12), (13) in~\cite{BoMo05}).
Using~\eqref{beta_n} and setting  $b_n:=\ln n$ we obtain that
\begin{align*}
\phi_{n,j} (\bt /\ln n)
 & =
 (1+o(1)) \exp\left\{  \frac{\alpha + 1+o(1)}{ \Si_1 (\xi)\ln n}
\sum_{r=0}^{n-1}\frac{\xi_r}{r}  (\bmu_r\bt^T
 + \eta_{r,n})  \right\}
 \\
 & =
  (1+o(1)) \exp\left\{  \frac{\alpha + 1+o(1)}{ \Si_1 (\xi) \ln n}
\left[  \sum_{r=0}^{n-1}\frac{\xi_r\bmu_r}{r} \, \bt^T
 +    \theta_n \sum_{r=0}^{n-1}\frac{\xi_r}{r} \right]\right\},
  \quad  |\theta_n |\le \eta_n.
\end{align*}

Next, applying to the sums in the exponential the following simple corollary of the
Stolz--Ces\`aro theorem: for a real sequence $\{x_n\},$ as $n\to\infty$,
\begin{align}
\mbox{if}\quad \frac1n \sum_{k=1}^n x_k \to x\in\R, \quad\mbox{then} \quad
 \frac1{\ln n} \sum_{k=1}^n \frac{x_k}k \to x
 \label{stolz}
\end{align}
(see, e.g., 2.3.25 in~\cite{KaNo00}), we obtain that $\lim_{(n,j)}\phi_{n,j} (\bt/\ln n)= \exp \{ i
\bla \bt^T\}$, where convergence is uniform in $\bt$ on any bounded subset of~$\R^d$, which
completes the proof of~\eqref{thm2_1}.

To prove~\eqref{thm2_2}, note that, for $\epsilon_n >0$ vanishing slowly enough (so that $U_{\epsilon_n n}/U_n = o(1)$ and $\ln\epsilon_n=o(\ln n)$; we will
assume without loss of generality that $\epsilon_n n\in\N$), setting $\bt_{r,n}:= (\bt \ln r)
/\ln n,$ one has from~\eqref{mean_m} that, for any fixed $\bt$, the ch.f.\ of $M_n
(\cdot \ln n)$ is equal to
\begin{align}
\varphi_n (\bt/\ln n)
%\frac1{U_n}\sum_{r=0}^n \sum_{j=1}^{k_r} u_{r,j} \phi_{r,j} (\bt/\ln n) %  + o (1)
 & = \frac1{U_n}\sum_{r=\epsilon_n n}^n \sum_{j=1}^{k_r} u_{r,j} \phi_{r,j} (\bt_{r,n}/\ln r) + o(1)
 \notag
 \\
 & = \exp \{ i \bt  \bla^T\} + o(1), \quad n\to\infty,
 \label{chf_M}
\end{align}
since $\bt_{r,n}= \bt + o(1)$ uniformly in  $r\in [\epsilon_n n, n]$, and
in view of the above-mentioned uniformity of convergence of~$\phi_{n,j} (\bt/\ln n)$ to the
exponential function.

\smallskip

(ii)~Using the second order uniform integrability assumption on~$\mathscr{P},$ it is
not hard to verify that, as $\|\bs\|\to 0,$
\[
f_r (\bs)-1 = i\bmu_r\bs^T  - \frac12 \bs \bm_r\bs^T+ o (\|\bs\|^2)
\quad\mbox{uniformly in~$r\ge 1$.}
\]
Therefore, for any fixed $\bt\in\R^d$,  taking into account that
$\sup_r\|\bmu_r\|<c<\infty$ (due to uniform integrability),   as $n\to\infty$ and $b_n\to \infty,$
\begin{align}
\ln (1   +\pi_r (f_r (\bt/b_n) -1))
 & = \pi_r (f_r (\bt/b_n) -1)
 + O\left( \frac{\pi_r^2  }{b_n^2}\right)
  \notag
  \\
  &= \pi_r \biggl( \frac{i\bmu_r\bt^T}{b_n}
  -\frac{\bt \bm_r\bt^T}{2b_n^2} \biggr)
   + o\left(\frac{\pi_r  }{b_n^2}\right)
   \label{ln_1}
\end{align}
from~\eqref{beta_n}. Now, observing that $f_{n } (\bt/b_n)=1 + o(1)$ as
$n\to\infty$ (again owing to the uniform integrability assumption), and setting
$b_n:=\sqrt{\ln n},$ one can derive from the first line of~\eqref{forma_b} that $\phi_{n,j} \bigl( {\bt}/{\sqrt{\ln n}} \bigr)$ is equal to
\begin{equation}
\label{exp_fu}
 (1+o(1)) \exp\left\{ \frac{i \bt}{\sqrt{\ln n}} \sum_{r=0}^{n-1} \pi_r \bmu_r^T
 - \frac{1}{2 \ln n} \bt \left(\sum_{r=0}^{n-1}  \pi_r\bm_r\right)\bt^T
  + \frac{o(1)}{\ln n} \sum_{r=0}^{n-1}  \pi_r \right\}.
\end{equation}
The first term in the argument of the exponential function in~\eqref{exp_fu} is
\[
\frac{i \bt}{\sqrt{\ln n}} \sum_{r=0}^{n-1} \pi_r \bmu_r^T
  =
 {i \bt}\sqrt{\ln n}\times  \frac1{ \ln n } \sum_{r=0}^{n-1} \pi_r \bmu_r^T
  \equiv
  i \bnu_n \bt^T \sqrt{\ln n}.
\]
Again using the Ces\`aro summability assumptions, \eqref{beta_n} and \eqref{stolz}, it is not difficult
to show, similarly to our argument in the proof of part~(i),  that, as $n\to\infty,$
\[
\frac{1}{  \ln n}   \left(\sum_{r=0}^{n-1}  \pi_r\bm_r\right)  \to
 (\alpha+1)\dfrac{\Si_1 (\xi \bm)}{\Si_1 (\xi)} .
\]
Finally, the last term in the argument of the exponential function in~\eqref{exp_fu}  is negligibly
small compared to
\[
\frac{1}{\ln n} \sum_{r=0}^{n-1}   \pi_r \to \alpha+1,
\]
the last relation again following from  \eqref{beta_n} and
\eqref{stolz}. This establishes the desired convergence of the distributions of $(\bX_{n,j}-\bnu_n\ln n)/\sqrt{n}$.

Now turn to the asymptotic behavior of
$M_n (\, \cdot\,\sqrt{\ln n} +\bnu_n\ln n)$. Let $\epsilon_n >0$ be a sequence vanishing slowly enough (so that $U_{\epsilon_n n}/U_n = o(1)$ and $\ln\epsilon_n=o(\sqrt{\ln n})$). Then, up to an additive term $o(1),$  the ch.f.\ $ \varphi_n (\bt /\sqrt{\ln n})e^{-i\nu_n\sbt^\top\sqrt{\ln n}} $ of that measure equals
\begin{align}
\frac1{U_n}& \sum_{r= \epsilon_n  n}^n \sum_{j=1}^{k_r} u_{r,j} \phi_{r,j} (\bt /\sqrt{\ln n})  e^{-i\sbnu_n\sbt^\top\sqrt{\ln n}}
\notag\\
 &=\frac1{U_n}\sum_{r=\epsilon_n  n}^n \sum_{j=1}^{k_r} u_{r,j} \phi_{r,j} (\bt_{r,n}/\sqrt{\ln r})  e^{-i\sbnu_r\sbt_{r,n}^\top\sqrt{\ln r}}
  \exp\biggl\{\frac1{\sqrt{\ln n}}\sum_{s=r}^{n-1}\pi_s \bmu_s \biggr\},
  \label{fi_M}
\end{align}
where $\bt_{r,n}:= \bt\sqrt{\frac{\ln r}{\ln n}}\to \bt $ uniformly in $r\in [\ep n, n]$ as $n\to\infty$ and, as $\sup_n\|\bmu_n\|\le c<\infty$ due to the uniform integrability assumption, one has from \eqref{beta_n} that, for some $c_1<\infty,$
\[
\biggl\|\sum_{s=r}^{n-1}\pi_s \bmu_s\bigg\|
 \le c \sum_{s=\epsilon_n n}^{n-1}\pi_s
 = c \frac{\al+1}{\Si_1 (\xi)} \sum_{s=\epsilon_n n}^{n-1}\frac{\xi_s}{s}(1+o(1))
 \le c_1    \sum_{s=\epsilon_n n}^{n-1}\frac{\xi_s}{s}.
\]
Setting $\Xi_n:=\frac1n \sum_{k=1}^{n } \xi_k $, note that $\xi_s = s\Xi_s- (s-1) \Xi_{s-1}$ and $\Xi_n\to\Si_1 (\xi)$ as $n\to\infty,$ and so
\begin{align*}
\sum_{s=\epsilon_n n}^{n-1}\frac{\xi_s}{s}
 & = \sum_{s=\epsilon_n n}^{n-1}\biggl(\Xi_s- \frac{s-1}s \, \Xi_{s-1}\biggr)
= \Xi_{n-1} - \Xi_{\epsilon_n n-1}+\sum_{s=\epsilon_n n}^{n-1} \frac1s\,  \Xi_{s-1}
\\
& = O\biggl(\sum_{s=\epsilon_n n}^{n-1} \frac1s\biggr)= O\bigl(|\ln \epsilon_n|\bigr) = o \bigl(\sqrt{\ln n}\bigr)
\end{align*}
by assumption. Therefore the last factor on the right-hand side of \eqref{fi_M} tends to one, so that the previously established convergence of the distributions of  $(\bX_{n,j}-\bnu_n\ln n)/\sqrt{n}$
completes the proof of Theorem~\ref{Thm2}.
\hfill$\Box$

\medskip

{\em Proof of Theorem~\ref{Thm3}.} (i)~All the computations below will   be valid on the event $A= \{\langle \widetilde\xi_1\rangle >0, \langle \tau_1\rangle >0\}.$

First set for brevity $a:=\langle \tau_1\rangle$ and
prove that all $\widetilde\zeta_n$ are finite on $A.$  By stationarity,   it suffices to show that  $\widetilde\zeta_0<\infty$ on~$A$. Using  Markov's inequality:
\[
\pr (\widetilde\xi_{-m} >e^{am/2}|\Ii)\le \langle\widetilde\xi_1\rangle e^{-am/2}, \quad m\ge 0,
\]
and the Birkhoff--Khinchin  theorem implying that $S_{0,m}=am (1+o(1))$ a.s.\ as $m\to\infty,$ we see from Borel--Cantelli's lemma that
\[
\widetilde\zeta_0 = \sum_{m=0}^\infty \widetilde\xi_{-m}e^{am(1+o(1))}
 \le \sum_{m=0}^\infty  e^{am/2+o(m)} <\infty \qquad \mbox{a.s.\ on $A$.}
\]

Setting $\vartheta_n:=\xi_n/\widetilde\xi_n -1$ (which is $o(1)$ a.s.\ by assumption \eqref{cond_xi}), we have
\begin{align}
W_n  &
 = \sum_{r=0}^n \xi_r e^{S_r} = e^{S_n}\sum_{m=0}^n \xi_{n-m}e^{-S_{n,m}}
 \notag\\
 & =   e^{S_n}\biggl(\widetilde\zeta_n + \sum_{m=0}^n \vartheta_{n-m} \widetilde \xi_{n-m} e^{-S_{n,m}} - e^{-\tau_0- S_n}\widetilde\zeta_{-1}\biggr).
 \label{Wn}
\end{align}
Note that here $e^{-S_n}= o(\widetilde\zeta_n)$ a.s. Indeed,  again by the Birkhoff--Khinchin theorem,
\[
\max_{m\le n/2} S_{n,m} =  \max_{m\le n/2} (an - a(n-m)) + o(n)=an/2 + o(n),
\]
and so
\begin{align}
\widetilde\zeta_n
 & \ge \sum_{m=0}^{n/2}\widetilde\xi_{n-m} e^{-S_{n,m}} \ge e^{- an/2 + o(n)} \sum_{m=0}^{n/2}\widetilde\xi_{n-m}
 \notag\\
  & = e^{- an/2 + o(n)}  \frac{n}2 \langle\widetilde\xi_1\rangle\gg  e^{- an  + o(n)}=  e^{-S_n} .
 \label{tildeze}
\end{align}
Moreover,
\[
\sum_{m=0}^n \vartheta_{n-m} \widetilde \xi_{n-m} e^{-S_{n,m}}
 = \sum_{0\le  m < 2n/3} + \sum_{ 2n/3\le m \le n}  =:\Sigma_1 + \Sigma_2,
\]
where $|\Sigma_1|\le \max_{r>n/3}|\vartheta_r|\widetilde\zeta_{n}=o(\widetilde\zeta_{n})$ a.s. Using notation $\overline \vartheta :=\sup_{j\ge 0} |\vartheta_{j}| $ (note that $\overline \vartheta< \infty$ a.s.) and the observation that $\min_{2n/3\le m \le n} S_{n,m}= 2na/3 + o(n)$ a.s., we obtain
\[
|\Sigma_2 |\le  \overline \vartheta\sum_{ 2n/3\le m \le n}\widetilde \xi_{n-m} e^{-S_{n,m}}
 \le \overline \vartheta  e^{-2na/3 + o(n)}  \sum_{ r=0}^{n/3}\widetilde \xi_{r} =  \overline \vartheta  e^{-2na/3 + o(n)}  \frac{n}3 \langle\widetilde\xi_1\rangle = o(\widetilde\zeta_n)
\]
from \eqref{tildeze}. Using the above bounds in~\eqref{Wn} we see that $W_n= \widetilde\zeta_n e^{S_n} (1+o(1))$ a.s.\ and hence, as $n\to\infty$,
\begin{equation}
\label{beta_n_3}
\pi_n = \frac{\widetilde\xi_n}{\widetilde\zeta_n}\,(1+o(1))\quad \mbox{a.s.}
\end{equation}

Now the conditional version of \eqref{forma_b} with $b_n:=n$ together with the uniform integrability assumption yields that, as $n\to\infty,$
\[
\phi_{n,j|\Vi} (\bt /n)=(1+o(1)) \exp\left\{  \frac{1+o(1)}{n} \sum_{r=1}^{n-1} \frac{\widetilde\xi_r}{\widetilde\zeta_r} (1+\eta_{r,n}^{(1)})(\widetilde\bmu_r\bt^T
 + \eta_{r,n}^{(2)})  \right\},
\]
where  $\max_{r\le n}|\eta_{r,n}^{(j)}|\le \eta_n=o(1)$ a.s., $j=1,2$. The desired convergence \eqref{Pn} follows now from   the Birkhoff--Khinchin theorem.

The assertion concerning the convergence of $M_{n|\Vi} (\,\cdot\, n)$ follows from the observation that the ch.f.\ of that measure, similarly to~\eqref{chf_M}, is equal to
\[
  \frac1{U_n}\sum_{r=0}^n \sum_{j=1}^{k_r} u_{r,j} \phi_{r,j|\Vi} \biggl(\frac{\bt}{r}\,\frac{r}{n}  \biggr)
  = \frac1{U_n}\sum_{r=0}^n   \sum_{j=1}^{k_r} u_{r,j} \exp \{ i \bt  \bla^T r/n+ \eta_{r,n}\},
\]
where $\eta_{r,n}=o(1)$  as $r\to\infty$, due to convergence~\eqref{Pn}.

\medskip

(ii) We start as in the proof of Theorem~\ref{Thm2}(ii), but, since   $\pi_n\not\to 0$ now, instead of~\eqref{ln_1} one has to use, for $b_n\to\infty,$ the following expansion:
 \begin{align*}
\ln (1    & +\pi_r (f_{r } (\bt/b_n)    -1))
  = \pi_r (f_{r } (\bt/b_n) -1)
 - \frac{\pi_r^2}2 (f_{r } (\bt/b_n) -1)^2 (1+o(1))
  \notag
  \\
   &= \pi_r \biggl( \frac{i\bmu_{r|\Vi}\bt^T}{b_n}
  -\frac{\bt \bm_{r|\Vi}\bt^T}{2b_n^2}  +o(b_n^{-2}) \biggr)
  - \frac{\pi_r^2(i\bmu_{r|\Vi} \bt^T +o(1))^2 }{2b_n^2}
    (1+o(1))
  \\
  & =\frac{i}{b_n} \pi_r \bmu_{r|\Vi}\bt^T
   -  \frac{ \pi_r}{2b_n^2}\bt \bigl(\bm_{r|\Vi} - \pi_r \bmu_{r|\Vi}^T \bmu_{r|\Vi}+ o(1) \bigr)\bt^T (1+o(1)),
\end{align*}
where the $o(1)$-terms are uniform in $r$ due to the uniform integrability assumptions.  Letting $b_n:=\sqrt{n},$ we obtain (cf.~\eqref{exp_fu}) that $\phi_{n|\Vi}
(\bt/\sqrt{n})$ is  given by
\begin{equation*}
  \exp\left\{ \frac{i \bt}{\sqrt{ n}} \sum_{r=0}^{n-1} \pi_r \bmu_{r|\Vi}^T
 - \frac{1}{2  n} \bt \left(\sum_{r=0}^{n-1} ( \pi_r\bm_{r|\Vi}-\pi_r^2 \bmu_{r|\Vi}^T\bmu_{r|\Vi})\right)\bt^T
  + \frac{o(1)}{\ln n} \sum_{r=0}^{n-1}  \pi_r +o(1)\right\}.
\end{equation*}
Therefore, using \eqref{beta_n_3}, we see that the conditional (given $\Vi$) ch.f.\ of $(\bX_{n,j} - \bka_n)/\sqrt{  n}$ is equal to
\[
\exp\left\{
 - \frac{1}{2  n} \bt
  \sum_{r=0}^{n-1} \frac{\widetilde \xi_r}{\widetilde \zeta_r} (1 + \eta^{(1)}_{r,n} )
  \Bigl(\widetilde \bm_r - \frac{\widetilde \xi_r}{\widetilde \zeta_r}\widetilde\bmu_{r }^T\widetilde\bmu_{r } +  {\beeta}^{(3)}_{r,n}  \Bigr) \bt^T
  + \frac{o(1)}{  n} \sum_{r=0}^{n-1}   \frac{\widetilde \xi_r}{\widetilde \zeta_r} (1+o(1)) +o(1)\right\},
\]
where the terms  $\eta^{(1)}_{r,n}$ and ${\beeta}^{(3)}_{r,n}$ are small uniformly  in~$r\le n$. Now the assertion of part~(ii) follows from the Birkhoff--Khinchin theorem.
\hfill$\Box$

%\medskip
%The author is grateful to the referees whose thoughtful comments helped him to improve the exposition of the paper.

\end{document}